\documentclass[a4,11pt]{article}%
\usepackage{graphicx}
\usepackage{amsmath}
\usepackage{amsfonts}
\usepackage{amssymb}
\usepackage{enumerate}
\usepackage{dsfont}
\usepackage{subfigure}

\newtheorem{theorem}{Theorem}

\newtheorem{claim}[theorem]{Claim}

\newtheorem{corollary}{Corollary}

\newtheorem{definition}{Definition}
\newtheorem{example}{Example}

\newtheorem{lemma}{Lemma}

\newtheorem{proposition}{Proposition}
\newtheorem{remark}{Remark}

\newcommand{\ep}{\varepsilon}
\newcommand{\dN}{\mathbb{N}}
\newcommand{\dR}{\mathbb{R}}

\newcommand{\rmd}{{\rm{d}}}

\newcommand{\cone}{{\rm{Cone}}}

\newcommand{\calC}{{\cal C}}

\newcommand{\calS}{\mathcal{S}}

\newcommand{\eps}{\varepsilon}

\newcommand{\conv}{{\rm conv}}
\newcommand{\sig}{{\sigma}}

\renewcommand{\phi }{\varphi }

\newcommand{\Cone}{\ensuremath{\rm Cone}}
\newcommand{\Conv}{\ensuremath{\rm Conv}}

\newcommand{\numbercellong}[2]
{
\begin{picture}(40,20)(0,0)
\put(0,0){\framebox(40,20)} \put(20,10){\makebox(0,0){#1}}
\put(33,15){\makebox(0,0){#2}}
\end{picture}
}

\newcommand{\numbercellongg}[1]
{
\begin{picture}(60,20)(0,0)
\put(0,0){\framebox(60,20)} \put(30,10){\makebox(0,0){#1}}
\end{picture}
}

\newcommand{\ignore}[1]{}
\def\qed{\vrule height 1.ex width 0.9ex depth -.1ex}
\newenvironment{proof}[1][Proof]{\textbf{#1.} }{\ \hfill  \qed}
\newenvironment{proofclaim}[1][Proof of Claim]{\textbf{#1.} }{\ \hfill \qed}
\usepackage{color}

\usepackage{graphicx}
\usepackage{epsfig}

\begin{document}
\title{Attainability in Repeated Games with Vector Payoffs%
\thanks{This research was supported in part by the Google Inter-university center for Electronic Markets and
Auctions. Lehrer acknowledges the support of the Israel Science Foundation, Grant \#538/11.
Solan acknowledges the support of the Israel Science Foundation, Grant \#212/09. Venel acknowledges the support of the Israel Science Foundation, Grant \#1517/11.}}
\author{Dario Bauso\thanks{Dipartimento di Ingegneria Chimica, Gestionale, Informatica, Meccanica, Universit\`a di Palermo, viale delle Scienze 90128, Palermo, Italy. e-mail: dario.bauso@unipa.it}~, Ehud Lehrer\thanks{The School of
Mathematical Sciences, Tel Aviv University, Tel Aviv 69978, Israel
and INSEAD, Boulevard de Constance, 77305 Fontainebleau, France. e-mail: lehrer@post.tau.ac.il}\ ,
Eilon Solan\thanks{The School of Mathematical Sciences, Tel Aviv
University, Tel Aviv 69978, Israel. e-mail: eilons@post.tau.ac.il}\ ,  and Xavier Venel \thanks{The School of Mathematical Sciences, Tel Aviv
University, Tel Aviv 69978, Israel. e-mail: xavier.venel@gmail.com}}
\maketitle

\begin{abstract}
We introduce the concept of attainable sets of payoffs in two-player
repeated games with vector payoffs. A set of payoff vectors is
called {\em attainable} by a player if there is a finite horizon $T$ such that the player can guarantee that after time $T$
the distance between the set and the cumulative
payoff is arbitrarily small, regardless of the strategy Player 2 is using.
We provide a necessary and sufficient condition for the attainability of a convex set, using the concept of $B$-sets.
We then particularize the condition to the case in which the set is a singleton,
and provide some equivalent conditions.
We finally characterize when all vectors are attainable.
\end{abstract}

\noindent\textbf{Keywords:} Attainability, continuous time, repeated games,
vector-payoffs, dynamic games, approachability.

\noindent {\emph{JEL} classification: C73, C72}

\thispagestyle{empty}

\bigskip

\newpage

\section{Introduction}

In various dynamic situations the stage-payoff is
multidimensional, and the goal of the decision maker is to drive
the total vector-payoff as close as possible to a given target
set. One such example is dynamic network models, which include a
variety of logistic applications such as production, distribution
and transportation networks. In the literature on dynamic network
flow control \cite{BBP10,BBP06,BMU00,BRU97,KT04}, the supplier
tries to meet a multidimensional demand. His goal is to ensure
that the difference between the {\em total} demand and the {\em
total} supply converges with time to a desirable target. One can
model such a situation as a two-player repeated game, where Player
1 is the decision maker and Player 2 represents the adversarial
market that controls demand.  In the distribution
network scenario,  for instance, the supplier has a desirable multidimensional
inventory level that he would like to maintain, despite erratic
behavior of the demand side. Having to deal with an adversarial
opponent requires the supplier to cope with the worst possible scenario.
This motivates our main objective: to find conditions that
characterize when a specific target set can be attained under
any possible demand pattern exhibited by the market.

A second example is the Capital Adequacy Ratio. The third Basel
Accord states that (a) the bank's Common Equity Tier 1 must be at
least 4.5\% of its risk-weighted assets at all times, (b) the bank's
Tier 1 Capital must be at least 6.0\% of its risk-weighted assets at
all times, and (c) the total capital, that is, Tier 1 Capital plus
Tier 2 Capital, must be at least 8.0\% of the bank's risk weighted
assets at all times. To accommodate this example in our setup,
consider the following 3-dimensional vector. The first coordinate
stands for the per-period difference between the bank's Common Equity Tier 1 and 4.5\%
of its risk-weighted assets; the second coordinate stands for
the per-period difference between the bank's Tier 1 Capital and 6.0\% of its
risk-weighted assets; and the third stands for the per-period difference
between the total capital and 8.0\% of the bank's risk weighted
assets. According to the Capital Adequacy Ratio the coordinates of
this vector should be nonnegative. Here, Player 1 represents the
bank's managers who control its assets, and Player 2 represents
market behavior, which is unpredictable and thought of as
adversarial. Thus, the goal of Player 1 is to design a strategy that
would drive the 3-dimensional total payoff to the target set -- the
nonnegative orthant. To ensure that they fulfill the requirements
of the Basel Accord, banks try to hold a capital buffer on top of
the regulatory minimum, and they periodically adjust their assets to
be at the top of the buffer \cite{HPS, vanRoy}.

To model such situations we study two-player repeated games with vector-payoffs in continuous time.
We say that a set $A$ in the payoff space is {\em
attainable} by Player 1 if there is a time $T$ such that for every level of proximity, $\ep>0$, Player 1 has a strategy guaranteeing that against
every possible strategy for Player 2,
the distance between $A$ and the cumulative payoff up to any time $t$ greater than $T$ is smaller than $\ep$.
If a set $A$ is attainable by a supplier, then in order to ensure that the inventory level would converge to $A$,  he can plan his actions based on historical inventory levels and market data.

The definition of attainability is close in spirit to the concept of
approachable sets \cite{B56}, which refers to the average stage-payoff rather than to the cumulative one.
While a set $A$ is attainable by Player 1 if he can ensure that
the cumulative payoff converges to it,
it is approachable by him if he can ensure that the average payoff converges to the set.

To illustrate the difference between these two notions, suppose that the target set of a supplier consists of one point, say $x$. If this set is attainable by him, it implies that the long-run inventory level is stable around $x$. On the other hand, if it is approachable, it merely guarantees that the average inventory level
converges to $x$.
This may happen also when the actual level itself does not converge to $x$,
and even when any fixed running average does not converge to $x$. This observation suggests that although
the notions of attainability and approachability are close to each other,  the flavor of the results and
their proofs are completely different.

One of our main results  characterizes attainable convex sets. It uses the concept of $B$-sets (see \cite{B56}).
It states that a convex set $Y$ is attainable by a player if and only if there exist two $B$-sets $C$ and $C'$ for that player
(or, alternatively, two approachable sets) and a nonnegative real number $\alpha$ such that
$\alpha C + \cone(C') \subseteq Y$.
The idea behind this result lies on two main properties that the cumulative payoff along
any possible trajectory of the game must have. The first is that the cumulative payoff must reach the set $Y$ within a certain time $T$, independently of the strategy of the other player.
The second property is that it has to remain close to $Y$ at any time after $T$.
The first property is responsible for the existence of $C$ and its role while the second property for that of $C'$.

In the case where $Y$ is compact then $\cone(C')$ is necessarily compact, which may happen only if $C'$ consists only of $\vec 0$. Our characterization entails that $\{\vec 0\}$ must then be approachable.
This observation enables us to provide necessary and sufficient conditions for a singleton (i.e., a set containing a single vector) to be attainable by Player~1.
We use it to show that {\em every} singleton is attainable by Player~1 if and only if the value of all scalar games
obtained from the vector-payoff game by projecting the payoff function on any direction $\lambda \neq \vec 0$ is positive.

%
%

\bigskip

The results presented here apply to games played in continuous time and where
players are allowed to use a special type of behavior strategies. These strategies are
characterized by an increasing sequence of positive real numbers that divide the time span $[0,\infty)$
into subintervals.
The play of a player in each interval depends on the play of the other player
{\em before} this interval starts and is independent on the other player's play during this interval.
This is equivalent to saying that before the game starts, a player sets an alarm clock to ring
at certain pre-specified times, and whenever the clock rings, the player
looks at the historical play path and determines how to play until the next time the clock rings.
In the literature of differential games this type of strategies is called
nonanticipating strategies with delay. We later discuss the interpretation of this type of strategies.

There is a
literature on decision problems related to dynamic multiinventory
in continuous time (see for instance, the continuous-time control strategy in \cite{BMU00}).
The control literature up to this point refers to one-person (the controller) decision problems with uncertainty.
To the best of our knowledge, this paper is the first that takes a strategic approach to these problems.

The paper is organized as follows. In Section \ref{sec:motivations} we provide
a motivating example. In Section \ref{sec:attainability} we introduce the model and main definitions.
In Section \ref{sec: results} we present our results,
and Section \ref{sec:comments} is devoted to discussing a few aspects related to the definition of
attainability and to the type of strategies that we are using.
Proofs are relegated to Section \ref{sec:proofs}.

\section{A motivating example}\label{sec:motivations}

\begin{figure}[ht]
\centering
\subfigure[Three distribution flows $f_A$, $f_T$, $f_B$
and two vendors requests $w_A$, $w_B$.]
{
\includegraphics[scale=.8]{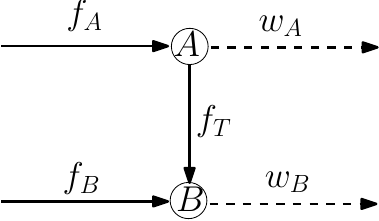}
\label{fig:subfig1}
}
\hskip 1.3cm
\subfigure[Factory manager can sell directly to vendors: node $C$ represents factory.]
{
\includegraphics[scale=.8]{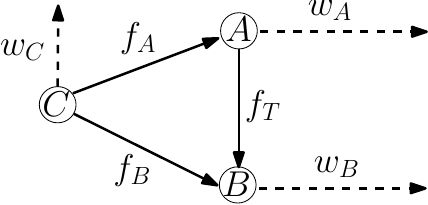}
\label{fig:subfig2}
}
\caption{Distribution network with warehouses $A$ and $B$.}
\end{figure}

This section details one motivation of our study: distribution networks.
Consider a distributer of a certain product who has two warehouses $A$ and $B$ in different regions.
Every month the distributer can order products from factories to each of the warehouses,
and he can transport products between the two warehouses,
while vendors order products from the warehouses.
This situation is described graphically in Figure \ref{fig:subfig1}.



 In Figure \ref{fig:subfig1}, $f_A$ and $f_B$ are the number of
products that are sent from factories to the two warehouses $A$ and
$B$, $f_T$ is the number of products that are transported from
warehouse $A$ to warehouse $B$, and $w_A$ and $w_B$ are the number
of products sent from the two warehouses to vendors. Negative flows
are interpreted as flows in the opposite direction; e.g., if vendors
return products to warehouse $A$ (resp. to warehouse $B$), then $w_A$ (resp. $w_B$) is negative.
If products are transported from warehouse $B$ to warehouse $A$,
then $f_T$ is negative. We analyze this situation in continuous
time. The change of stock in the two warehouses is given by the
2-dimensional vector

\begin{eqnarray*}
 u(a_1^t,a_2^t) =
 \underbrace{\left(\begin{array}{ccc} 1 & -1 & 0 \\
 0 & 1 & 1\end{array}\right)}_{F}
 \underbrace{\left(\begin{array}{c} f_A^t \\f_T^t\\f_B^t\end{array}\right)}_{a_1^t} -
 \underbrace{\left(\begin{array}{c} w_A^t\\w_B^t\end{array}\right)}_{a_2^t},
 \end{eqnarray*}
where $a_1^t = (f_A^t,f_B^t,f_T^t)$ is the decision variable of the distributer,
and $a_2^t = (w_A^t,w_B^t)$ is the uncontrolled market demand at time $t$.

Suppose that the number of products that can be ordered by vendors  at each time instance
is bounded by 2, and the number of products that can be returned by vendors to each warehouse at every time instance is 3. In other words,  $w_A^t$ and $ w_B^t$ are in $[-3,2]$.
Suppose also that the amount of product that the distributer can order from or return to the factories and
transport between the two warehouses is bounded by $5$.

This situation can be described by a two-person game as follows.
The distributer (Player 1) has 8 actions
\[ (5,5,5), (5,5,-5), (5,-5,5), (5,-5,-5), (-5,5,5), (-5,5,-5), (-5,-5,5), (-5,-5,-5), \] while
the market demand or nature (Player 2) has 4 actions
\[ (-3,-3), (-3,2), (2,-3), (2,2). \] The payoffs correspond to the change of stock in the two warehouses, and are given by the following table:


\centerline{
\begin{picture}(280,190)(-40,0)
\put(-40,8){(-5,-5,-5)}
\put(-40,28){(-5,-5,5)}
\put(-40,48){(-5,5,-5)}
\put(-40,68){(-5,5,5)}
\put(-40,88){(5,-5,-5)}
\put(-40,108){(5,-5,5)}
\put(-40,128){(5,5,-5)}
\put(-40,148){(5,5,5)}
\put(20,170){(-3,-3)}
\put(80,170){(-3,2)}
\put(140,170){(2,-3)}
\put(200,170){(2,2)}
\put( 0,0){\numbercellongg{(3,-7)}{}}
\put( 0,20){\numbercellongg{(3,3)}{}}
\put( 0,40){\numbercellongg{(-7,3)}{}}
\put( 0,60){\numbercellongg{(-7,13)}{}}
\put( 0,80){\numbercellongg{(13,-7)}{}}
\put( 0,100){\numbercellongg{(13,3)}{}}
\put( 0,120){\numbercellongg{(3,3)}{}}
\put( 0,140){\numbercellongg{(3,13)}{}}
\put( 60,0){\numbercellongg{(3,-12)}{}}
\put( 60,20){\numbercellongg{(3,-2)}{}}
\put( 60,40){\numbercellongg{(-7,-2)}{}}
\put( 60,60){\numbercellongg{(-7,8)}{}}
\put( 60,80){\numbercellongg{(13,-12)}{}}
\put( 60,100){\numbercellongg{(13,-2)}{}}
\put( 60,120){\numbercellongg{(3,-2)}{}}
\put( 60,140){\numbercellongg{(3,8)}{}}
\put(120,0){\numbercellongg{(-2,-7)}{}}
\put(120,20){\numbercellongg{(-2,3)}{}}
\put(120,40){\numbercellongg{(-12,3)}{}}
\put(120,60){\numbercellongg{(-12,-13)}{}}
\put(120,80){\numbercellongg{(8,-7)}{}}
\put(120,100){\numbercellongg{(8,3)}{}}
\put(120,120){\numbercellongg{(-2,3)}{}}
\put(120,140){\numbercellongg{(-2,13)}{}}
\put(180,0){\numbercellongg{(-2,-12)}{}}
\put(180,20){\numbercellongg{(-2,-2)}{}}
\put(180,40){\numbercellongg{(-12,-2)}{}}
\put(180,60){\numbercellongg{(-12,8)}{}}
\put(180,80){\numbercellongg{(8,-12)}{}}
\put(180,100){\numbercellongg{(8,-2)}{}}
\put(180,120){\numbercellongg{(-2,-2)}{}}
\put(180,140){\numbercellongg{(-2,8)}{}}
\end{picture}
}

\centerline{Figure 2: The strategic-form game corresponding to the situation.}

\bigskip
\noindent
At every time instance the two players choose their actions.
Each market behavior translates into a mixed action of Player 2,
and each behavior of the distributer corresponds to a mixed action of Player 1.
The (2-dimensional) total payoff up to time $t$ is the number of products that are stored in each of the two warehouses.
The goal of the distributer is to ensure that the total number of products in each warehouse does not exceed its capacity,
that is, that the total payoff should not exceed a certain (2-dimensional) bound.


Figure \ref{fig:subfig2} describes the case where the factory manager can sell directly to vendors,
bypassing the distribution to warehouses. This situation can be represented by adding an additional node $C$
modeling the factory, and an edge that represents the market demand.
The stock is now a 3-dimensional vector, as we have to take into account the inventory available at the
factory, and consequently the change in the stock modifies as shown below:
\begin{eqnarray*}
 u(a_1^t,a_2^t) =
\left(\begin{array}{ccc} 1 & -1 & 0 \\
 0 & 1 & 1\\ -1 & 0 & -1 \end{array}\right)
\left(\begin{array}{c} f_A^t \\f_T^t\\f_B^t\end{array}\right) -
 \left(\begin{array}{c} w_A^t\\w_B^t \\ w_C^t\end{array}\right).
 \end{eqnarray*}

A recurrent question in the network flow control literature \cite{BBP10,BBP06,BMU00,BRU97,KT04}
is  about conditions that ensure the existence of a control strategy
that drives the excess supply vector to a desired target level in $\mathbb R^m$
regardless of the unpredictable realization of the demand.
The equivalence between the excess supply and the cumulative payoff in the dynamic game motivates our study. The rest
of the paper is devoted to the analysis of conditions under which Player 1
has a strategy ensuring the attainability of a convex set,
regardless of the behavior of Player 2.

Situations where the target is to control the total payoff occur also in
production and transportation networks. Production networks describe
production processes and activities necessary to turn raw materials
into intermediate products and eventually into final products. The
nodes of the networks represent raw materials and intermediate/final
products. The buffer at each single node $i$ models the amount of
material or product of type $i$ stored or produced up to the current
time, and hyper-arcs describe the materials or products consumed
(tail nodes) and produced (head nodes) in each activity or process.
Transportation networks model the flow of commodities, information,
or traffic; nodes of the networks represent hubs and the buffers at
the nodes describe the quantity of commodities present in the hubs.
The edges describe transportation routes.

\subsection{Related control and optimization literature}
We highlight two main streams of related literature, one from the control area and the second from the optimization area.
These two bodies of literature have two main elements in common: i) the interest towards robustness, and ii) the presence of
a network dynamic flow scenario.

Connections between robust control and noncooperative game theory has a long history (see, e.g., \cite{BB91}).
Robust control is the area of control theory that looks for strategies that ``control''  the state of a
dynamical system, for instance, drive it to a given set,  despite the effects of disturbances (see the seminal paper \cite{BR71}).
Among the foundations of robust control we find two main notions that can be related
to attainability and are surveyed in \cite{B99}.
The first notion, \textit{robust global attractiveness},
refers to the property of a set to ``attract'' the state of the system under a proper control strategy,
independently of the effects of the disturbance.
The second notion, \textit{robustly controlled invariance},
describes the property of a set to bound the state trajectory
under a proper control strategy, independently of the effects of the disturbance.
Both notions are widely exploited in a variety of works that contribute to the use of robust control in dynamic network flow models \cite{BBP10,BBP06,BMU00,BRU97,KT04}.

A second stream of literature can be identified under the name of ``robust optimization''. This is a relatively recent technique that
describes uncertainty via sets and optimizes the worst-case cost over those sets (see, e.g., the introduction to the special issue
\cite{BGN06}). The use of robust optimization techniques in dynamic network models is the main focus of \cite{AP06,AZ,BT06}.
There, a main theme is to ``adjust'' some of the supplier's decision variables to the uncertain outcome.
More specifically, some variables are determined
before the outcome is realized while the rest are determined
after the outcome is realized.
Such a problem formulation is referred to as ``Adjustable Robust
Counterpart'' (ARC) problem, or ``two-stage robust optimization
with recourse'' and as it will be clear later
it shares striking similarities with the formulation of attainable strategies presented in the current paper.

This paper focuses on the game theoretic aspects related to
attainable sets. A discussion on applications to network flow control problems is
introduced in a companion paper \cite{LSB11}.

\section{Attainability}\label{sec:attainability}
In the first part of this section we introduce the mathematical model of repeated game in continuous time
and elaborate on the type of strategies used by the players.
In the remaining part, we provide a formal definition of attainability.

\subsection{The model}\label{sec:model} We study a two-player
repeated game with vector payoffs in continuous time $\Gamma$. The set of
players is $N = \{1,2\}$, and the finite set of actions of each
player $i$ is $A_i$. The instantaneous payoff is given by a function
$u: A_1 \times A_2 \to \dR^m$, where $m $ is a natural number. We
assume w.l.o.g.\ that payoffs are bounded by 1, so that $u : A_1
\times A_2 \to [-1,1]^m$. We extend $u$ to the set of mixed action
pairs, $\Delta(A_1) \times \Delta(A_2)$, in a bilinear fashion.
The one-shot vector-payoff game $(A_1,A_2,u)$ is
denoted by $G$ and we will say that the game in continuous time $\Gamma$ is \emph{based on} $G$. For $i \in \{1,2\}$, $-i$ denotes the opponent of $i$.

The game $\Gamma$ is played over the time interval $[0,\infty)$.
We assume that the players use nonanticipating
behavior strategies with delay, which we define below. Roughly, a
nonanticipating behavior strategy with delay divides time
into intervals. The behavior of a player in a given interval depends on
the behavior of the other player up to the beginning of the interval.
In other words, the way a player plays during a given interval of time does
not affect the way the opponent plays during that interval. Still,
it may affect the other player's play in subsequent intervals.

Formally, denote by $\calC_i$ the set of all {\em controls} of player $i$,
that is, the set of all measurable functions from the
time space, $[0,\infty)$, to player $i$'s mixed actions. That is,
\[ \calC_i := \left\{ a_i : [0,\infty) \to \Delta(A_i), ~ a_i \hbox{ is measurable}\right\}. \]

\begin{definition}
\label{def strategy} A function $\sigma_i : \calC_{-i} \to \calC_i$ is a {\em behavior strategy
with delay} (or simply a \emph{strategy}) for player $i$,
if there
exists an increasing sequence of real numbers $(\tau_i^k)_{k \in \dN}$
such that for every $a_{-i}, a'_{-i} \in \calC_{-i}$,
\[ a_{-i}(t) = a'_{-i}(t) \ \ \ \forall t \in [0,\tau_i^k) \ \ \ \Longrightarrow
(\sigma_i(a_{-i}))(t) = (\sigma_i(a'_{-i}))(t) \ \ \ \forall k \in \dN, \forall t \in [\tau_i^{k},\tau_i^{k+1}),
\]
where $\tau_i^0 = 0$.
\end{definition}
In the sequel we refer to the real numbers $(\tau_i^k)_{k \in \dN}$
in Definition \ref{def strategy} as the \emph{updating times related to
$\sigma_i$}.

\begin{remark}
In the literature on differential games a strategy as the one defined above is called a \emph{\em nonanticipating strategy with delay}.
An equivalent formulation, which may look more transparent to game theorists, is as follows.
A strategy for player $i$ is a list $(\tau_i^k,\sigma_i^k)_{k \in \dN}$
where $(\tau_i^k)_{k \in \dN}$  is an increasing sequence of real numbers,
and for each $k \in \dN$,
$\sigma_i^k$ is a function that maps play paths (of both players) on the interval $[0,\tau_i^{k})$
to plays of player $i$ in the interval $[\tau^{k}_i,\tau^{k+1}_{i})$.

Later on when defining or referring to a strategy we use this formulation.
\end{remark}

\begin{remark}
Models with continuous time are typically used as a tractable version
of discrete-time models, where the gap between two consecutive
stages is small. This is the case here as well. Suppose that time
is discrete, and that the time interval between any two successive
decisions is extremely small. Suppose moreover, that observing
opponent's actions is time consuming and possibly costly. Thus,
players cannot observe each other's actions at every time. Rather, they
observe their opponent's actions relatively rarely compared to the
frequency in which actions are taken. Our continuous-time model captures
this aspect: although time is continuous,  players observe
historical play and decide how to play in the next interval in
discrete times. Neither updating nor new decision is taken place
between two updating times.
\end{remark}

Every pair of strategies $\sigma =
(\sigma_1,\sigma_2)$ uniquely determines a play path
$(a^\sigma(t))_{t \in \dR_+}$. The cumulative payoff-vector up to time $T$
associated with the pair of strategies $\sigma$ is given by
\begin{equation}
\label{equ96} \gamma^T(\sigma) = \int_0^T u(a^\sigma(t)) \rmd t \in
\dR^m.
\end{equation}
Sometimes we denote it by $\gamma^T_G$ when we wish to emphasize that the payoff is in the game based on $G$.
Note that since the payoffs are bounded by 1, the integral in (\ref{equ96}) is well-defined.

\subsection{Attainability: the definition}

The subject matter of this paper is the concept of attainable sets. A
set of vectors is attainable by a player if he can guarantee that the distance between the set and the cumulative
payoff converges to 0, regardless of the strategy of the opponent. We provide a definition here and two alternative
ones are discussed later in
Section \ref{sec:comments}.
\begin{definition}
\label{def attainability} A nonempty closed set  $Y \subseteq \dR^m$  is {\em attainable} by
Player 1 if there is $T > 0$ such that for every $\ep > 0$ there is
a strategy $\sigma_1$ of Player 1 such that\footnotemark
\[ d(\gamma^t(\sigma_1,\sigma_2),Y) \leq \ep, \ \ \ \forall t \geq T, \forall\sigma_2. \]
\end{definition}
\footnotetext{The distance referred to throughout the
paper is the Euclidean distance
and the norm is the $L_2$-norm
$\|.\|_2$. The distance between a point $x$ and a set $A$ is, therefore, $d(x,A)=\min_{y \in Y}
\|x-y\|_2$.}

A set $Y$ is attainable if there is a finite horizon $T$ such that Player 1 can
ensure, against any possible strategy of Player 2, that the cumulative
payoff up to any time $t\geq T$ is within $\ep$ from $Y$. Note that the
time $T$ is uniform across all levels of precision. That is, in order for $Y$ to be attainable by Player 1, Player
1 must be able to guarantee that the cumulative payoff at any time longer than
$T$ would be within any $\ep$ from $Y$. However, different $\ep$'s
might require different strategies employed by Player 1. It might
therefore happen that although $Y$ is attainable by Player 1, the cumulative payoff would never touch
$Y$ itself.
We say that the
strategy $\sigma_1$ in Definition \ref{def attainability}
\emph{attains the set $Y$ up to $\ep$}.

When $Y$ contains a single vector and it is attainable by Player 1, we say
that the vector $x$ is attainable by him. Denote by $W$ the set of attainable vectors.

Since the notion of attainability is related to that of approachability,
we recall the definition of the later (in a continuous-time framework). Denote the mean vector-payoff between times $0$ and $T$ by
\[
\overline{\gamma}^T(\sigma_1,\sigma_2)=\frac{1}{T} \gamma^T(\sigma_1,\sigma_2).
\]

\begin{definition}
A nonempty closed set $Y \subseteq \mathbb{R}^m$ is \emph{approachable by Player
$1$} if for every $\epsilon>0$ there exist $T>0$ and a strategy
$\sigma_1$ of Player $1$ such that
\[
d\left(\overline{\gamma}^t(\sigma_1,\sigma_2),Y \right) \leq
\epsilon, \ \ \ \forall t
\geq T, \forall \sigma_2.\]
\end{definition}
\ignore{
Differential games are typically defined by their dynamics. Here the dynamic of the average payoff is given by
\[
\dot{\overline{\gamma}}(t)=\frac{-\overline{\gamma}(t)+u(a^\sigma(t))}{t}.
\]
}
The original definition of approachable sets \cite{B56} was given in discrete time repeated games. A set is approachable in the discrete time model if and only if it is approachable in the continuous time one
 (see, \cite{SQS09}).

The definitions of attainability and  approachability are close in spirit. There is, however, a significant difference between the two concepts. A set is approachable if the {average}
payoff converges to it, while a set is attainable if the cumulative
payoff converges to the set. In other words, approachability refers
to the convergence of the \emph{average} payoff, while attainability
to the convergence of the \emph{cumulative} payoff.

\section{Results}\label{sec: results}

This section presents the main results. The first, Theorem 1, characterizes closed and convex attainable sets.
Using this result, we derive a characterization of attainable compact convex sets and of attainable vectors.
Finally, the last result provides a stronger condition that ensures that
any vector $x \in \dR^m$ is attainable.


%

For every nonempty closed set $Y \subseteq \mathbb{R}^m$ and every $z \in \mathbb{R}^m$ we denote by
$\Pi_Y(z):=\{y \in Y \mid d(z,y)=d(z,Y) \}$ the set of points in $Y$ closest to $z$.
When the set $Y$ is convex, $\Pi_Y(z)$ contains a single point.
Our main theorem characterizes closed convex attainable sets. To state the result we borrow from \cite{B56} the concept of $B$-set.
\begin{definition}
A nonempty closed set $C \subseteq \mathbb{R}^m$ is a \emph{$B$-set for Player~1} if for every
$z \in \mathbb{R}^m$ there exists $c\in \Pi_C(z)$ and $x \in \Delta(A_1)$ such
that
\[
\langle u(x,a_2)-c,z-c \rangle \leq 0, \ \ \  \forall a_2\in A_2.
\]
\end{definition}

If a nonempty closed set contains a $B$-set, then it is approachable \cite{B56}.
Conversely, every approachable set contains a $B$-set \cite{Hou,S02}.
Therefore, a set is approachable if and only if it contains a $B$-set.

We show that a closed convex set is attainable if and only if
it contains a certain sum of two $B$-sets.
For every set $Y \subseteq \mathbb{R}^m$ we denote the cone spanned by $Y$ by $\Cone(Y)=\{\alpha y,
\alpha \in \mathbb{R}^+, y\in Y\}$.

\begin{theorem}\label{maintheo}
A closed convex set $Y \subseteq \mathbb{R}^m$ is attainable by Player~1 if
and only if there exist $\alpha>0$ and two $B$-sets for that player $C$ and $C'$ such
that
\begin{equation}
\label{equ:main}
\alpha C + {\rm Cone}(C') \subseteq Y.
\end{equation}
\end{theorem}

The idea behind the theorem is that any trajectory attaining a set $Y$ consists of two parts.
A first sub-trajectory reaches $Y$ in a fixed finite time. This is represented by the $\alpha C$.
A second sub-trajectory stays close to $Y$.
When $Y$ is unbounded, keeping the trajectory within $Y$ is equivalent to keeping the direction in which the trajectory progresses within a proper range.
This is represented by the second term ${\rm Cone}(C')$. In the special case where $Y$ is compact (i.e.,  also bounded)
Player~1 can ensure that the trajectory will remain in $Y$ only if he can keep the trajectory close to $\vec 0$,
so that the set $\cone(\{\vec 0\})$ is attainable.
The inclusion in (\ref{equ:main})  is justified by the observtion that any superset of an attainable set is attainable too.

%
%

We deduce now several corollaries. In the first we focus on compact convex sets.
Whenever $Y$ is compact, the characterization of Theorem \ref{maintheo} can be simplified.
Indeed, since the only compact cone is $\{\vec 0\}$ whenever $Y$ is a compact convex attainable sets we must have $C = \{\vec 0\}$.
By Theorem \ref{maintheo} it follows that there exists $\alpha > 0$ and a $B$-set $C$ such that
$C \subseteq \frac{1}{\alpha}Y$.
By setting $\delta= \frac{1}{\alpha}$, we infer that $\delta Y$ is approachable.
Since every approachable set contains a $B$-set, this yields the following result.

\begin{corollary}
\label{coro1} A compact convex set $Y\subseteq \mathbb{R}^m$ is attainable by Player~1 if and only if
\begin{enumerate}
\item[\emph{\textbf{B1}}]   The vector $\vec 0 \in \mathbb{R}^m$ is approachable by Player~1, and
\item[\emph{\textbf{B2}}]   There exists a scalar $\delta>0$ such that $\delta Y$ is approachable by Player~1.
\end{enumerate}
\end{corollary}

The following example, borrowed from \cite{B56}, shows that the result introduced above does not hold when $Y$ is not convex.

\begin{example}
Consider the following payoff repeated game with 2-dimensional payoffs:

\centerline{
\begin{picture}(110,65)(-10,0)
\put(-10, 8){$B$}
\put(-10,28){$T$}
\put( 20,50){$L$}
\put(60,50){$R$}
\put( 0,0){\numbercellong{$(1,0)$}{}}
\put( 0,20){\numbercellong{$(0,0)$}{}}
\put( 40,0){\numbercellong{$(1,1)$}{}}
\put( 40,20){\numbercellong{$(0,0)$}{}}
\end{picture}
}
\centerline{Figure 3: The payoff function in Example \ref{sec: example}.}

\bigskip

Define $Y:=\{(1/2,t), \ t \leq 1/4\} \cup \{(1,t), \ t \geq 1/4\}$.
It was shown in Blackwell \cite{B56} that the set $Y$ is not approachable by Player~1.
One can verify that none of its dilatations is approachable by Player~1.
Nevertheless, the set $Y$ is attainable by Player~1. Indeed the following behavior strategy attains it for Player~1:
\[
\sigma_1(t) = \left\{
\begin{array}{lll}
B & \ \ \ \ \ & t \in [0,\frac{1}{2}),\\
T & & t \in [\frac{1}{2},1), \gamma_1^{1/2} \leq \frac{1}{4},\\
B & & t \in [\frac{1}{2},1), \gamma_1^{1/2} > \frac{1}{4},\\
T & & t \geq 1.
\end{array}
\right.
\]
\end{example}

In the particular case where $Y = \{x\}$ is a singleton, we can again be more precise. We seperate the attainability of $\vec 0$ and the attainability of $x \neq \vec 0$. To state the next result we need the following
notations. Let $\lambda\in \mathbb{R}^m$. Denote%
\footnote{ The inner product is defined by $\langle x,y\rangle := \sum_{i=1}^m x_iy_i$ for every $x,y \in \dR^m$.}
 by $\langle
\lambda, G\rangle$ the zero-sum one-shot game whose set of players and their
action sets are
as in the game $G$, and the payoff
that Player 2 pays to Player 1
is $\langle \lambda,
u(a_1,a_2) \rangle$ for every  $(a_1,a_2)\in A_1 \times A_2$. As a
zero-sum one-shot game, the game  $\langle \lambda, G\rangle$ has a value, denoted $v_{\lambda}$.

For every mixed action $p \in \Delta(A_1)$ denote
\[ D_1(p) = \{ u(p,q) \colon q \in \Delta(A_2)\}. \]
$D_1(p)$ is the set of all payoffs that might be realized when
Player 1 plays the mixed action $p$. If $v_{\lambda}\geq 0 $ (resp.
$v_{\lambda}>0 $), then there is a mixed action $p \in \Delta(A_1)$
such that $D_1(p)$ is a subset of the closed half space $\{x\in
\mathbb{R}^m\colon \langle \lambda, x\rangle \geq 0\}$ (resp. half
space $\{x\in \mathbb{R}^m\colon \langle \lambda, x\rangle > 0\}$).
Thus $D_1(p)$ and $\lambda$ are in the same half-space, or,
equivalently, $D_1(p)$ and $-\lambda$ are in two different
half-spaces.

\begin{corollary}
\label{coro2} The following three properties are equivalent.
\begin{enumerate}
\item[\emph{\textbf{C1}}]   The vector $\vec 0 \in \mathbb{R}^m$ is attainable by Player 1
\item[\emph{\textbf{C2}}]   The vector $\vec 0 \in \mathbb{R}^m$ is approachable by Player 1.
\item[\emph{\textbf{C3}}]   For every  $\lambda\in \mathbb{R}^m$, $v_\lambda \geq 0$.
\end{enumerate}
\end{corollary}

The equivalence between $C1$ and $C2$ is an immediate consequence of Corollary \ref{coro1}. Based on that $\vec 0$ is approachable if and only if it is a $B$-set, the equivalence between $C2$ and $C3$ follows from \cite{B56}.

%
%
%
%

The following result characterizes when a given vector $x \neq \vec 0$ is attainable. For every  $y\in\dR^m$
denote by $(G-y)$ the two-player one-shot game that is identical to $G$ except for
its payoff function. The payoff function of $(G-y)$ is $(u-y)$, where
$(u-y)(a_1,a_2)=u(a_1,a_2)-y$ for every $a_1 \in A_1$ and $a_2 \in A_2$.

\begin{corollary}\label{coro3} Let $\vec 0 \not =x \in \mathbb{R}^m$. The vector
$x$ is attainable by Player 1 if and only if
\begin{enumerate}
\item[\emph{\textbf{D1}}]   The vector $\vec 0 \in \mathbb{R}^m$ is attainable by Player 1
\end{enumerate}
and either one of the following conditions holds:
\begin{enumerate}
\item[\emph{\textbf{D2}}]   There is $\delta>0$ such that the vector $\vec 0 \in \mathbb{R}^m$ is attainable by Player 1
in the game based on $(G-\delta x)$.
\item[\emph{\textbf{D3}}]   There is $\delta>0$ such that, for every $\lambda\in \mathbb{R}^m$, $v_\lambda \geq \delta \langle x,\lambda\rangle$.
\item[\emph{\textbf{D4}}]   There is $\delta_0>0$ such that for every $q\in \Delta(A_2)$ there is $p\in \Delta(A_1)$
and $\delta > \delta_0$ satisfyingze
$u(p,q)=\delta x$.
\end{enumerate}
\end{corollary}

Conditions $\textbf{D2}$ and $\textbf{D3}$ are reformulations of $\textbf{B2}$. Condition $\textbf{D4}$ needs additional work in order to be proven. The proofs are differed to the last section.

Corollary \ref{coro3} implies that whenever any vector $x$ is attainable, so
is the vector $\vec 0$. Since attainability is concerned with the
cumulative payoff, once a target level is (almost) reached, this
level should be maintained in the long run. This means that once a
neighborhood of a target level $x$ is reached, from that point in time and on
$\vec 0$ ought to be attained. This is the reason why
$\vec 0$ is attainable when any vector $x$ is attainable, and why $\vec 0$ plays a
major role in the theory of attainability. However, Condition
\textbf{B1} alone is not sufficient for the attainability of other vectors other then $\vec 0$ itself.

\bigskip
The previous result  naturally leads to a sufficient condition for a vector to be attainable.
\begin{proposition}\label{suffvec}
Let $x \in \mathbb{R}^m$ such that
\begin{itemize}
  \item [\textbf{E1}] $v_{\lambda} \geq 0$ for every $\lambda \in \mathbb{R}^m \setminus \{\vec 0\}$, and
 \item [\textbf{E2}] For every $\lambda \in \mathbb{R}^m \setminus \{\vec 0\}$, if $\langle \lambda,x\rangle \geq 0$ then $v_{\lambda} > 0$.
   \end{itemize}
Then $x$ is attainable.
\end{proposition}

We deduce the following theorem which deals with the case where all the vectors are
attainable.
\begin{theorem}
\label{theorem2} The following statements are equivalent:
\begin{itemize}
  \item [\emph{\textbf{F1}}] Every vector $x \in \mathbb{R}^m$ is
    attainable by Player 1;
 \item [\emph{\textbf{F2}}] $v_{\lambda}>0 $ for
every  $\lambda\in \mathbb{R}^m \setminus \{\vec 0\}$.
   \end{itemize}
\end{theorem}

The fact that Condition $\textbf{F2}$ implies Condition $\textbf{F1}$ is a consequence of Proposition \ref{suffvec}.
Indeed given that $\textbf{F2}$ is true, then for every $x \in \mathbb{R}^m$, Condition $\textbf{E2}$ is satisfied and thus every $x \in \mathbb{R}^m$ is attainable.

The converse implication can be obtained by focusing on Condition $\textbf{D3}$ in Corollary~\ref{coro3}.
Assume that Condition \textbf{F1} holds.
For every $\lambda\in \mathbb{R}^m \setminus \{\vec 0\},$ the vector $x=\lambda$ is attainable and satisfies $\langle x,\lambda\rangle > 0$.
Therefore, Condition \textbf{D3} implies that $v_\lambda >0$, and Condition $\textbf{F2}$ holds as well.

\begin{remark}\label{rem delta}  \noindent
If Condition {\textbf{F2}} is satisfied, then for every open half
space $H$ of $\mathbb{R}^m$ there is a mixed action $p \in
\Delta(A_1)$ such that $D_1(p) \subseteq H$. Standard continuity and
compactness arguments imply that in this case there is $\delta_1 >
0$ such that for every half space $H$ there is $p \in \Delta(A_1)$
satisfying $d(D_1(p), H) \geq \delta_1$. Stated differently, there
is $\delta_2 > 0$ such that for every vector $\lambda$ whose
$\ell_1$-norm is $1$, $\langle \lambda,u(p,q) \rangle>\delta_2$ for
every $q\in \Delta(A_2)$.
\end{remark}

Note the difference between Condition {\textbf{C3}} of Corollary
\ref{coro2} and Condition {\textbf{F2}} of Theorem
\ref{theorem2}. In the former, the value of the scalar-payoff game
with payoffs $\langle \lambda,u(p,q) \rangle$ is nonnegative for
every direction $\lambda \in \dR^m \setminus \{\vec 0\}$, while in the latter it is
strictly positive. The former guarantees attainability of the vector $\vec 0$,
while the latter guarantees that every vector is attainable.

\section{Discussion}\label{sec:comments}
The model and the results described above give rise to a number of
additional questions. (a)  What are the analogous results in  discrete time repeated game to the ones we obtained? (b)
Are there different notions of attainability that do not impose a uniform time of convergence? (c)
What happens if the updating times are not predetermined and can be selected as a function of the information
available up to the updating time?

We next elaborate on these questions and highlight a few
open problems left for future research.

\subsection{Continuous time versus discrete time.}
The characterization presented in Theorem \ref{maintheo}
depends crucially on the continuous time setting. The following example shows
that it is invalid when time is discrete.

\begin{example} \normalfont
\label{example1}
Consider a game in \emph{discrete} time where payoffs are
one-dimensional and each player has two actions. Payoffs are given
by the following matrix:

\centerline{
\begin{picture}(110,65)(-10,0)
\put(-10, 8){$B$} \put(-10,28){$U$} \put( 20,50){$L$} \put(
60,50){$R$}
\put( 0,0){\numbercellong{$2-1$}{}}
\put( 0,20){\numbercellong{$-2-1$}{}}
\put( 40,0){\numbercellong{$2+1$}{}}
\put( 40,20){\numbercellong{$-2+1$}{}}
\end{picture}
\begin{picture}(130,65)(-40,0)
\put(-40,18){$=$}
\put(-10, 8){$B$}
\put(-10,28){$U$}
\put( 20,50){$L$}
\put( 60,50){$R$}
\put( 0,0){\numbercellong{$1$}{}}
\put( 0,20){\numbercellong{$-3$}{}}
\put( 40,0){\numbercellong{$3$}{}}
\put( 40,20){\numbercellong{$-1$}{}}
\end{picture}
}
\centerline{Figure 3: The payoff function in Example \ref{example1}.}

\bigskip

The payoffs in this game are the sum of two numbers, one determined by
Player 1 (-2 if he plays $U$, 2 if he plays $B$), and the other by
Player 2 (-1 if she plays $L$, 1 if she plays $R$).

Condition \textbf{C3} is satisfied, and therefore $0$  is attainable
by Player 1. The following strategy guarantees that the cumulative
payoff is
within $9 \cdot
2\eta$ from $0$ at any $t>2$, where $\eta > 0$ is given; the details of the proof can be found
in the proof of Theorem \ref{maintheo}. Divide the time line into
countably many blocks, where the length of the $k$-th block is
$\eta\over k$. In the $k$-th block Player 1 plays $U$ if the
cumulative payoff at the beginning of the block is positive, and he
plays $B$ otherwise.

We show that $0$ is not attainable by Player 1 in the game in discrete time.
When time is discrete, a behavior strategy for a player is a function that assigns
a mixed action to each past history.
For every $\ell \in \dN$, let $p^\ell$ be the mixed action played by Player 1 at stage $\ell$.
Note that $p^\ell$
depends on past play.
Let $\sigma_2$ be the strategy that at each stage $\ell$ plays $L$ if
$p^\ell(U) \geq \frac{1}{2}$, and $R$ otherwise.
The stage payoff is then at least 2 whenever Player 2 plays $R$,
and at most $-2$ whenever Player 2 plays $L$.
In particular, if the total payoff up to stage $\ell$ is in the interval $[-\frac{1}{2},\frac{1}{2}]$,
then the payoff up to stage $\ell+1$ lies outside this interval. Thus, the cumulative
 payoff does not converge to 0.\hfill \qed
\end{example}

Example \ref{example1} suggests that the characterization of
the set of attainable vectors in games in discrete time
is more challenging than the characterization in continuous time.

\subsection{Alternative definitions of attainability}

We here provide two alternative definitions of the concept of attainability,
which we term asymptotic attainability and weak asymptotic attainability.
We then explore some relations between the three definitions.

For every set $Y \subseteq \dR^m$ we denote by $B(Y,{\eps})$ the set
of all points whose distance from at least one point in $Y$ is less
than $\eps$, that is,
\[ B(Y,\ep) := \{ x\in \dR^m \colon d(x,Y) < \ep\}. \]
When $Y$ contains a single point $x$, we write $B(x,\ep)$ instead of
$B(\{x\},\ep)$.

\begin{definition} \label{def weak}
(i) The set $Y \subseteq \dR^m$ is {\em asymptotically attainable}
by Player 1 if there is a strategy $\sigma_1$ for Player 1 such that
for every strategy $\sigma_2$ of Player 2,
\begin{equation}
\label{equ def attainability}
 \lim_{T \to \infty} d(\gamma^T(\sigma_1,\sigma_2),Y) = 0.
\end{equation}\\
(ii) The set  $Y$ is {\em weakly asymptotically attainable} by
Player 1, if the set $B(Y,{\eps})$ is asymptotically attainable by Player
1 for every $\eps>0$.
\end{definition}

Asymptotic attainability requires that a set is asymptotically
reached by the cumulative payoff without putting any bound on the time it takes to reach the set.
Attainability, on the other hand, requires that a set is
approximately reached in a bounded time, independent of the degree of
approximation. Weak asymptotic attainability relaxes both time boundedness and the level of the approximation precision.
A set $Y$ is weakly asymptotically
attainable if any neighborhood $B(Y,{\eps})$ of $Y$ can be
asymptotically attained, without having a universal bound on the
time at which this neighborhood is reached.

Any attainable set is also weakly asymptotically attainable and any asymptotically attainable set is weakly asymptotically attainable as well.
In addition, observe that the set of asymptotically attainable vectors
and the set of weakly asymptotically attainable vectors are convex cones.
The definition implies that the set of weakly attainable vectors is also closed.

Using Corollary \ref{coro3} , we now show that attainability of a vector does not imply its asymptotic attainability.
This implies in particular that these two concepts are not identical.

\begin{example}\label{sec: example} \normalfont
We provide an example where the vector $\vec 0$ is attainable but
not asymptotically attainable. Consider the following
game where payoffs are 2-dimensional, each player has 2 actions, and the payoffs are scalar and given by:

\centerline{
\begin{picture}(110,65)(-10,0)
\put(-10, 8){$B$}
\put(-10,28){$U$}
\put( 20,50){$L$}
\put(60,50){$R$}
\put( 0,0){\numbercellong{$0$}{}}
\put( 0,20){\numbercellong{$1$}{}}
\put( 40,0){\numbercellong{$-1$}{}}
\put( 40,20){\numbercellong{$0$}{}}
\end{picture}
}
\centerline{Figure 4: The payoff function in Example \ref{sec: example}.}

\bigskip

In this game $v_{\lambda}=0$ for every $\lambda \in \dR$. Thus,
for every $\lambda \in \dR^2$ one has $v_{\lambda}\ge 0$, and therefore  corollary \ref{coro3} implies that the vector $\vec 0$ is attainable by
Player 1. We argue that $\vec 0$ is not asymptotically attainable by
Player 1. Assume that Player 1 implements a strategy $\sigma_1$. In an initial time
interval the strategy $\sigma_1$ plays one of the rows with a
positive probability. Consider the strategy $\sigma_2$ of Player 2
that plays constantly a column that generates a nonzero vector in
that initial interval. For instance, if $\sigma_1$ plays the action $U$ with
positive probability in the initial time interval, then $\sigma_2$
play the action $L$ always. The initial period produces a
nonzero payoff and this payoff is not diminishing to zero because
Player 2 keeps playing the same column forever. This example shows
that $\vec 0$ is attainable by Player~1 but not asymptotically attainable by
him.

We point out that the argument mentioned above shows in fact that
$\vec 0$ is not attainable in the corresponding game in discrete
time as well.
\end{example}

The following example shows that a weakly attainable vector need not be attainable.
\begin{example}  \normalfont
\label{example4}
Consider a two-player game where payoffs are 2-dimensional,
Player 1 has 3 actions, Player 2 has 2 actions, and the payoff
function is given by the left-hand side matrix in Figure 5.

\centerline{
\begin{picture}(90,110)(-10,-20)
\put(-10, 8){$B$}
\put(-10,28){$M$}
\put(-10,48){$U$}
\put( 20,70){$L$}
\put( 60,70){$R$}
\put( 0,0){\numbercellong{$(0,0)$}{}}
\put( 0,20){\numbercellong{$(0,0)$}{}}
\put( 0,40){\numbercellong{$(1,1)$}{}}
\put( 40,00){\numbercellong{$(0,0)$}{}}
\put( 40,20){\numbercellong{$(1,1)$}{}}
\put( 40,40){\numbercellong{$(0,1)$}{}}
\put(11,-20){The game $G$}
\end{picture}
}

\centerline{Figure 5: The payoff function of the game $G$ in Example \ref{example4}.}

\bigskip

The vector $(0,0)$ is attainable by Player 1, using the strategy
that always plays $B$. The vector $x:=(1,1)$ is weakly asymptotically
attainable according to Definition \ref{def weak}. Indeed, given
$\ep > 0$ consider the strategy $\sigma_1^\ep$, with updating times
$(\tau_1^k)_{k \in \dN}$ defined by $\tau_1^k = k\ep$ for $k \in \dN$, that is
defined as follows.
\begin{itemize}
\item   If the total payoff up to time $\tau_1^k$ is not in the set $B((1,1),\ep)$,
during the time interval $[\tau_1^k,\tau_1^{k+1})$ play the mixed action $[\ep(U),(1-\ep)(M)]$.
\item   If the total payoff up to time $\tau_1^k$ is in the set $B((1,1),\ep)$,
during the time interval $[\tau_1^k,\tau_1^{k+1})$ play the action $B$.
\end{itemize}
For every $t \geq \frac{1}{\ep}$ one has $d(\gamma^t(\sigma^\ep_1,\sigma_2),(1,1)) < \ep$,
so that the vector $x$ is indeed weakly asymptotically attainable by Player 1.

The vector $x$, however, is not attainable by Player 1 (according to
Definition \ref{def attainability}).
To show this claim we use Corollary \ref{coro3}
and prove that Condition \textbf{D4} does not hold for $x$.
Indeed, fix $\delta_0 > 0$, and set $q := [(1-\frac{\delta_0}{2})(L),\frac{\delta_0}{2}(R)]$.
Let $p \in \Delta(A_1)$ be arbitrary.
If $u(p,q) = \delta x = (\delta,\delta)$ for $\delta > 0$, then necessarily $p_{_U} = 0$.
One can verify that $u(p,q)$ cannot be equal to $\delta x$ for $\delta > \delta_0$,
and therefore Condition \textbf{D4} does not hold for $x$.
\end{example}

\begin{remark}
The proof of Theorem \ref{theorem2} shows that every vector $x \in \dR^m$ is
attainable by Player 1 if and only if every vector $x \in \dR^m$ is asymptotically attainable by Player 1.
Example \ref{sec: example} shows that attainability does not imply
asymptotic attainability. We are unable to tell whether or not
asymptotic attainability implies attainability.
\end{remark}

\subsection{Alternative strategies in continuous time.}

The strategies we use here are nonanticipating strategies with delay. In these
strategies the times $(\tau_i^k)_{k \in \dN}$ at which a player observes past play
are independent of the play of the other player.
One could consider a broader class of strategies in which
$(\tau_i^k)_{k \in \dN}$ are stopping times. In other words,
$\tau_i^{k+1}$ is a time that depends on (that is, it is measurable with
respect to) the information available to player $i$ at time $\tau_i^k$, for each
$k \in \dN$. In this type of strategies, the updating times $(\tau_i^k)_{k \in \dN}$,
are not predetermined real numbers, as in Definition \ref{def strategy}.
Our results remain valid even if Player 2 is allowed to use a
strategy from this broader class of strategies.

 \subsection{Additional open problems}
%
%

The results above refer to attainability of a convex set, and did not discuss attainability, asymptotic attainability, or weak asymptotic attainability of nonconvex sets. Characterizing when a set of payoffs is attainable (according to these three definitions) remains open. We also leave attainability in discrete time and attainability when payoffs are discounted for future investigations.

\section{Proofs}\label{sec:proofs}

\subsection{Proof of Theorem \ref{maintheo}}

The aim of this section is to prove the characterization of attainable closed convex sets.
We first provide the outline of the proof.

\subsubsection{Outline}
Continuous approachability and discrete approachability are equivalent \cite{SQS09}.
This justifies the use of the notion of $B$-sets in the study of games in continuous time.

Given $\alpha>0$, a closed convex set $Y \subseteq \dR^m$, and two $B$-sets $C,C' \subseteq \dR^m$, two $B$-sets such
that
\[ \alpha C +
{\rm Cone}(C') \subset Y,
\]
we check that, by convexity of $Y$, we can replace the $B$-sets with
their convex hull which are approachable:
\[
\alpha \Conv(C) + \Cone(\Conv(C')) \subset Y.
\]

Then we prove that given two approachable convex sets $C$ and $C'$
the set $\alpha C + {\rm Cone}(C')$ is attainable. Thus $Y$ is
attainable.

To show the converse implication, given a set $Y$, we define the
family of sets $\overline{Y}_t$ as the intersection of $[-1,1]^m$
and $\frac{1}{t}Y.$ We show that the family $\overline{Y}_t$ admits
limit values $\overline{Y}_\infty.$ Moreover, for every $t>0,$ we
have
\[
t \overline{Y}_t + \Cone(\overline{Y}_\infty) \subset Y
\]

We prove that there exists $T \in \dR_+$ such that for all $t\geq T$,
the set $\overline{Y}_t$ is approachable. This implies that $\overline{Y}_T$
and $\overline{Y}_\infty$ are approachable and each one contains a
$B$-set \cite{S02}.
It follows that (\ref{equ:main}) holds with $C=\overline{Y}_T$ and $C' = \overline{Y}_\infty$.

\subsubsection{The condition is sufficient}

Let $Y$ be a closed convex set.
Suppose that there exists $\alpha>0$ and two $B$-sets
$C$ and $C'$ such that
\[
\alpha C + {\rm Cone}(C') \subset Y.
\]
We prove that $Y$ is attainable.
Since any superset of an approachable set is approachable,
the sets $\Conv(C)$ and $\Conv(C')$ are approachable.
Since $Y$ is convex, these two sets are subsets of $Y$ and satisfy
\[
\alpha \Conv(C) + \Cone(\Conv(C')) = \Conv (\alpha C + \Cone(C'))
\subset \Conv(Y) = Y.
\]

We now prove the following.

\begin{proposition}\label{toatt}
Let $\alpha>0$ and $C,C'$ be two closed convex approachable subsets of $\mathbb{R}^m$. Then
$\alpha C+{\rm Cone}(C')$ is attainable.
\end{proposition}

The following lemma claims that the distance between a point and a set is a convex function.
For every finite collection $(C_i)_{i=1}^n$ of nonempty subsets of $\dR^m$
and every collection $(\lambda_i)_{i=1}^m$ of scalars,
denote
\[ \sum_{i=1}^n \lambda_i C_i:=\{z \in \mathbb R^m \mid z=\sum_{i=1}^n  \lambda_i c_i, \, \forall c_i\in C_i,\, \forall i=1,\ldots,n\}. \]

\begin{lemma}
Let $n \in \dN$,
let $(x_i)_{i=1}^n$ be points in $\dR^m$,
and let $(C_i)_{i=1}^n$ be nonempty closed subsets of $\dR^m$.
For every collection of positive real numbers $(\lambda_i)_{i=1}^n$ one has
\[
d\left( \sum_{i=1}^n \lambda_i x_i, \sum_{i=1}^n \lambda_i C_i \right) \leq \sum_{i=1}^n \lambda_i d(x_i,C_i).
\]
\end{lemma}

\begin{proof}
For every $i \in \{1,2,\ldots,n\}$ let $c_i$ be a point in $C_i$ that satisfies $d(x_i,c_i) = d(x_i,C_i)$.
Note that $c := \sum_{i=1}^n \lambda_i x_i \in \sum_{i=1}^n \lambda_i C_i$.

\begin{eqnarray*}
d\left( \sum_{i=1}^n \lambda_i x_i, \sum_{i=1}^n \lambda_i C_i \right) &\leq& \left\|  \sum_{i=1}^n \lambda_i x_i - c\right\| \\
                      &\leq& \sum_{i=1}^n \lambda_i \|x_i-c_i \| \\
                      &=& \sum_{i=1}^n \lambda_i d(x_i,C_i).
\end{eqnarray*}
\end{proof}

When $C$ is convex, $\sum_{i=1}^n \lambda_i C  = \left(\sum_{i=1}^n \lambda_i\right) C$,
and therefore we obtain the following corollary.

\begin{corollary}
\label{cor3}
Let $n \in \dN$,
let $(x_i)_{i=1}^n$ be points in $\dR^m$,
and let $C$ be a nonempty closed and convex subset of $\dR^m$.
For every collection of positive real numbers $(\lambda_i)_{i=1}^n$ one has
\[
d\left( \sum_{i=1}^n \lambda_i x_i, \left(\sum_{i=1}^n \lambda_i\right) C \right) \leq \sum_{i=1}^n \lambda_i d(x_i,C).
\]
\end{corollary}

For every strategy $\sigma_i$ of Player~$i$ let $\sigma^{\beta}_i$ be the strategy $\sigma_i$
accelerated by a factor $\beta$. That is,
$(\sigma^{\beta}_i(a_{-i}))(t) := (\sigma_i(\widehat a_{-i}))(\beta t)$,
where $\widehat a_{-i}(t) = a_{-i}(\beta t)$.

The following result, which holds since time is continuous,
states that if the strategy $\sigma_1$ of Player~1 guarantees that the distance between the average payoff up to time $t$ and a given set $C$ is less than $\ep$,
then the accelerated strategy $\sigma^\beta_1$ ensures that the distance between the average payoff up to time $\frac{t}{\beta}$ and the set $\frac{C}{\beta}$
is at most $\frac{\epsilon}{\beta}$. 

\begin{lemma}
\label{lemma:1}
Let $C \subseteq \dR^m$, $T > 0$, and $\ep > 0$.
If the strategy $\sigma_1$ satisfies
\[
d\left({\overline{\gamma}^t(\sigma_1,\sigma_2)},C\right)
\leq \epsilon, \ \ \ \forall t \geq T, \ \forall \sigma_2,
\]
then the strategy $\sigma^\beta_1$ satisfies
\[
d\left({\overline{\gamma}^t(\sigma^\beta_1,\sigma_2)},\frac{C}{\beta}\right)
\leq \frac{\epsilon}{\beta}, \ \ \ \forall t \geq \frac{T}{\beta}, \ \forall \sigma_2.
\]
\end{lemma}

\begin{proof}
For every strategy $\sigma_2$ of Player 2 one has
\begin{eqnarray}
\nonumber
\gamma^t(\sigma^{\beta}_1,\sig_2)
&=&
\int_0^t u(a^s(\sigma_1^\beta,\sigma_2))\rmd s\\
\label{equ97}
&=&
\frac{1}{\beta}\int_0^{\beta t} u(a^s(\sigma_1,\sigma_2^{1/\beta}))\rmd s\\
&=&\frac{1}{\beta}\gamma^{\beta t}(\sigma_1,\sig^{1/\beta}_2).
\nonumber
\end{eqnarray}
We deduce that
\[ d\left({\overline{\gamma}^t(\sigma_1^\beta,\sig_2)},\frac{x}{\beta}\right) \leq \frac{\eps}{\beta}, \ \ \ \forall t \geq \frac{T}{\beta}, \forall \sigma_2, \]
as desired.
\end{proof}

A corollary of this result, which has its own interest but will not be used here,
is that the set of attainable vectors is a convex set.

A second corollary of Lemma \ref{lemma:1} is that if $C$ is a set that is approachable by Player~1,
then for every $t > 0$ he can ensure that the total payoff up to time $t$ is arbitrarily close to $tC$.

\begin{corollary}
\label{corollary:5}
Let $C$ be a set that is approachable by Player~1. For every
$\delta>0$ and every $s \in (0,1)$ there exists a strategy $\sigma^*_1$ such that
\[
d\left(\gamma^s(\sigma^*_1,\sigma_2),s C \right) \leq \delta, \ \ \ \forall \sigma_2.
\]
\end{corollary}

\begin{proof}
Because the set $C$ is approachable by Player~1,
there is a constant $K$, a strategy $\sigma_1$ and $T > 0$ such that
\[
d\left(\overline\gamma^t(\sigma_1,\sigma_2),C \right) \leq \frac{K}{\sqrt{t}}, \ \ \ \forall \sigma_2, \forall t \geq T.
\]
W.l.o.g.~we can assume that $T \geq \left(\frac{K}{\delta}\right)^2$.
It follows that
\[
d\left(\gamma^t(\sigma_1,\sigma_2),tC \right) \leq t\frac{K}{\sqrt{t}}, \ \ \ \forall \sigma_2, \forall t \geq T,
\]
and by Lemma \ref{lemma:1} we have for every $\beta > 0$,
\[
d\left(\gamma^{\beta t}(\sigma^\beta_1,\sigma_2),\beta tC \right) \leq \beta t\frac{K}{\sqrt{t}}, \ \ \ \forall \sigma_2, \forall t \geq T.
\]
For every  $s \in (0,1)$ substitute $t := T$ and $\beta := \frac{s}{T}$, to obtain
\[
d\left(\gamma^{s}(\sigma^\beta_1,\sigma_2),sC \right) \leq s\frac{K}{\sqrt{T}}, \ \ \ \forall \sigma_2.
\]
The result follows since $s \leq 1$ and $T \geq \left(\frac{K}{\delta}\right)^2$.
\end{proof}

We now strengthen Corollary \ref{corollary:5} to prove that if a set $C$ is approachable by Player~1,
then he can guarantee that the total payoff remains close to the cone generated by it, for every $s \geq 0$.
\begin{lemma}
\label{lemma:5}
Let $C$ be a closed and convex set that is approachable by Player~1. For every
$\varepsilon>0$ there exists a strategy $\sigma^*_1$ such that
\[
d\left(\gamma^s(\sigma^*_1,\sigma_2),s C \right) \leq \epsilon, \ \ \ \forall \sigma_2, \forall s \geq 0.
\]
\end{lemma}

\begin{proof}
The strategy $\sigma_1$ is given by concatenating strategies that satisfy Corollary \ref{corollary:5},
with properly chosen $s$'s and $\ep$'s.

Let us start by fixing $\ep>0$.
For each $k \in \dN$ let $\sigma_1^k$ be a strategy that satisfies Corollary \ref{corollary:5} with $s=\frac{\ep}{k}$
and $\delta = \frac{\ep}{2^k}$.
Set $t_0 := 0$ and
\[ t_{k+1} := t_k + \frac{\ep}{k}, \ \ \ \forall k \geq 0. \]
Let $\sigma_1$ be the strategy of Player~1 that for each $k \geq 0$,
at time $t_k$ forgets past play and follows the strategy $\sigma_1^k$ until time $t_{k+1}$.

By construction and by Corollary \ref{cor3}, for each $k \geq 0$ we have
\[ d(\gamma^{t_k}(\sigma_1,\sigma_2),t_k C) \leq \sum_{j=1}^k \frac{\ep}{2^k} \leq \frac{\ep}{2}. \]
Since $t_{k+1}-t_k \leq \frac{\ep}{2}$, and since payoffs are bounded by 1,
the triangle inequality implies that
\[ d(\gamma^{t}(\sigma_1,\sigma_2),t C) \leq \ep, \ \ \ \forall t \geq 0, \]
as desired.
\end{proof}

\bigskip

We are now ready to complete the proof of Proposition  \ref{toatt}.

\bigskip

\begin{proof} [Proof of Proposition \ref{toatt}] Let $\epsilon>0$, we
build a strategy $\sigma_1^*$ such that
\[
d\left(\gamma^t(\sigma^*_1,\sigma_2),\alpha C+ \Cone(C') \right)
\leq 2\epsilon, \ \ \ \forall t\geq \alpha,\forall \sigma_2.
\]

We define $\sigma_1(C)$ and $\sigma_1(C')$   -- two strategies given by
the previous lemma applied respectively to $C$ and $C'$. We
define a strategy $\sigma^*$ of Player $1$:  follow $\sigma_1(C)$ until time
$\alpha$ and then follow $\sigma_1(C')$.

Let $t\geq \alpha$ and $\sigma_2$ a strategy of Player $2$. By
construction, we have
\begin{align*}
d\left(\gamma^t(\sigma^*_1,\sigma_2),\alpha C+ (t-\alpha)(C')
\right) & \leq d\left(\gamma^\alpha(\sigma^*_1(C),\sigma_2),\alpha C
\right)+
d\left(\gamma^{t-\alpha}(\sigma^*_1(C'),\sigma'_2),(t-\alpha)(C')
\right)\\
 & \leq 2\epsilon,
\end{align*}
where $\sigma'_2$ being the continuation strategy of Player $2$ after time
$\alpha$. The set $\alpha C + \Cone(C')$ is attainable.
\end{proof}

\subsubsection{The condition is necessary}

In this section we prove that if a closed and convex set $Y$ is
attainable by Player~1, then there exists $\alpha>0$ and two $B$-sets $C$ and $C'$ such
that
\[
\alpha C + {\rm Cone}(C') \subset Y.
\]

For every $t \geq 0$ define
\[ Y_t := Y \cap [-t,t]^m \]
and
\[ \overline{Y}_t:=\frac{1}{t}Y_t \subseteq [-1,1]^m. \]

The Hausdorff metric is a metric over closed subsets of $\dR^m$, and defined as follows:
\[ d^H(X,Y) := \sup_{x \in X} \inf_{y \in Y} \|x-y\|. \]
It is well known that the set of closed subsets of a compact set is compact in this metric.
 This implies that the sequence
$(\overline{Y}_t)_{t > 0}$ has an accumulation point, with respect to the
Hausdorff metric.
\begin{lemma}\label{coneinc}
Let $Y$ be a nonempty closed convex set and  $\overline{Y}_\infty$ be an accumulation point of the sequence
$(\overline{Y}_t)_{t > 0}$.
If  $\vec 0 \in Y$ then
\[{\rm
Cone}(\overline{Y}_\infty) \subseteq Y.\]
\end{lemma}

%
%

\begin{proof}
For each $t \geq 0$ the set $\overline{Y}_t$ is a compact subset of $[-1,1]^m$.
The first claim follows by observing that the collection of compact subsets of $[-1,1]^m$ is itself compact in the Hausdorff metric.

We now turn to the second claim.
To show that ${\rm Cone}(\overline{Y}_\infty) \subset Y$ we fix a point $z \in
{\rm Cone}(\overline{Y}_\infty)$ and construct a sequence of points $(x'_{n})_{n\in \dN}$ in $Y$ that converges to it.
Since $Y$ is closed, this will prove that $z \in Y$.
Because $z \in \cone(\overline{Y}_\infty)$,
there exist $\alpha>0$ and $y\in
\overline{Y}_\infty$ such that $z= \alpha y$.

Since the sequence $(\overline{Y}_t)_{t > 0}$ converges to $\overline Y_\infty$ in the Hausdorff metric,
for every $n \in \dN$ there exists $t_n \geq \alpha$ satisfying
\[
d^H(\overline{Y}_{t_n},\overline{Y}_\infty) \leq \frac{1}{n}.
\]
In particular, there is $w_n \in \overline{Y}_{t_n}$ such that $d(w_n,y) \leq \frac{1}{n}$.
Setting $x_n := t_nw_n \in Y_{t_n}$ we deduce that $d\left(\frac{1}{t_n}x_n,y\right) \leq \frac{1}{n}$,
or equivalently,
\[ d\left(\frac{\alpha}{t_n}x_n,\alpha y\right) \leq \frac{\alpha}{n}. \]
Since (i) $\alpha \leq t_{n}$, (ii) $x_n,\vec 0 \in Y$, and (iii) $Y$ is convex,
it follows that $x'_n := \frac{\alpha}{t_n}x_n$ is in $Y$, and the result follows.
\end{proof}

\begin{lemma}\label{decomp}
Let $Y$ be a nonempty closed and convex set,
and let $\overline{Y}_\infty$ be an (Hausdorff metric) accumulation point of the sequence $(\overline{Y}_t)_{t > 0}$.
For every $t\geq 1$ one has
\[
t \overline{Y}_t + {\rm Cone}(\overline{Y}_\infty) \subseteq Y.
\]
\end{lemma}

\begin{proof}
For every $y \in \dR^m$ denote by $Z^y=Y-y$.
Then $\lim_{t \to \infty} d(\overline{Z^y}_t,\overline{Y}_t) = 0$,
so that $\overline{Y}_\infty$ is an accumulation point of the sequence $(\overline{Z^y}_t)_{t > 0}$.

Fix now $y \in Y$.
The set $Z^y$ is nonempty, closed, convex, and contains $\vec 0$,
so that by Lemma \ref{coneinc}
\[ \cone(\overline{Y}_\infty) = \cone(\overline{Z^y}_\infty) \subseteq Z^y = Y-y. \]
In particular,
\[ y + \cone(\overline{Y}_\infty) \subseteq Y. \]
The result  follows from the fact that this inclusion holds for every $y \in Y$,
and because $t \overline{Y}_t \subseteq Y$ for every $t > 0$.
\end{proof}

To conclude the proof that the condition is necessary we show that the sets $\overline{Y}_t$ are
approachable by Player 1, provided $t$ is large enough.
This will imply that the set
$\overline{Y}_\infty$, as an accumulation point of approachable sets, is approachable itself.

We start by finding a condition, lightly weaker than that in the definition of approachable sets, which is equivalent to it.

\begin{lemma}
\label{lemma:7}
A nonempty closed set $Y$ is approachable by Player~1 if and only if for every $\epsilon>0$ there exists a strategy $\sigma_1$
and $T > 0$ such that
\begin{align*}\label{eq1}
d \left(\overline{\gamma}^T(\sigma_1,\sigma_2),Y
\right) \leq \epsilon, \ \ \ \forall \sigma_2.
\end{align*}
\end{lemma}

\begin{proof}
The fact that if $Y$ is approachable by Player~1 then it satisfies the condition in the lemma follows from the definition of approachability.
For the converse implication, fix $\ep > 0$, and let $\sigma_1$ and $T$ be the strategy of Player~1 and the positive real number that are given
by the condition in the lemma.
Let $\sigma'_1$ be the strategy of Player~1 that plays in blocks of length $T$;
at the beginning of each block the strategy forgets past play and starts implementing $\sigma_1$ anew.
The reader can verify that $\sigma'_1$ approaches $Y$.
\end{proof}

We are finally ready to prove that the condition in Theorem \ref{maintheo} is necessary.
Let $Y$ be a  closed and convex set that is attainable by Player~1.
Therefore, there exists $T > 0$ such
that for every $\epsilon>0$ there exists a strategy $\sigma_1$ satisfying
\[
d(\gamma^t(\sigma_1,\sigma_2),Y)\leq \epsilon, \ \ \ \forall t\geq T, \ \forall \sigma_2.
\]
By Lemma \ref{lemma:7} this implies that the set $\frac{1}{t}Y$ is approachable by Player~1, provided that $t \geq T$.
Since payoffs are bounded by 1, it follows that the set $\overline{Y}_t = \frac{1}{t}Y \cap [-1,1]^m$ is also approachable by Player~1, provided that $t \geq T$.
Finally, the definition of approachability implies that the set $\overline{Y}_\infty$,
as the Hausdorff limit of sets which are approachable by Player~1, is also approachable by Player~1.
Since every set that is approachable by Player~1 contains a $B$-set for that player (\cite{Hou,S02}),
the proof of the necessity of the condition is complete.

\subsection{Proof of Corollary \ref{coro3}}

By Corollary \ref{coro1}, the set $\{x \}$ is attainable by Player~1 if and only if $\vec 0$ is approachable by Player~1 and
there exists $\delta>0$ such that the vector $\delta x $ is approachable by Player~1.
We will show that the second property is equivalent to the Conditions $\textbf{D2}$ and $\textbf{D3}$.
We will then prove that, given that $\textbf{D1}$ is satisfied, we can replace any one of these conditions with $\textbf{D4}$.

\noindent\underline{Part 1:}
$\delta x$ is approachable by Player~1 if and only if Condition \textbf{D2} holds.

Note that the vector $\delta x $ is approachable by Player~1 in the game with matrix payoff $G$
if and only if the vector $\vec 0$ is approachable by him in the game with matrix payoff $G-\delta x$.
Since $\vec 0$ is approachable by a player if and only if it is  attainable by him, the result follows.

\noindent\underline{Part 2:}
$\delta x$ is approachable by Player~1 if and only if Condition \textbf{D3} holds.

Let us write the $B$-set condition with respect to the singleton $\delta x$.
The vector $\delta x$ is approachable by Player~1 if and only if it is a $B$-set for him, that is,
\begin{align*}
\forall z\in \mathbb{R}^m,\ \exists x\in\Delta(A_1)\  \forall y\in \Delta(A_2) \langle u(x,y)-\delta x, z-\delta x\rangle \leq 0.
\end{align*}

Setting $\lambda=\delta x-z$, we obtain
\begin{align*}
\forall \lambda \in \mathbb{R}^m,\ \exists x\in\Delta(A_1)\  \forall y \in \Delta(A_2)&  \langle u(x,y), \lambda \rangle \geq \langle\delta x, \lambda\rangle,\\
& \langle u(x,y), \lambda \rangle \geq \delta \langle x, \lambda\rangle,
\end{align*}
which is equivalent to
$v_\lambda \geq \delta \langle x, \lambda\rangle$ for every $\lambda \in \mathbb{R}^m$,
which is Condition
$\textbf{D3}$

\noindent\underline{Part 3:} If the vector $x$ is attainable by
Player 1, then Condition \textbf{D4} is satisfied.

Suppose to the contrary that Condition \textbf{D4} is not satisfied.
That is, for every $\delta_0 > 0$ there is $q \in \Delta(A_2)$ such
that for every $p \in \Delta(A_1)$ one has $u(p,q) \neq \delta x$
for every $\delta > \delta_0$. We divide the argument into two
cases.

\bigskip
\noindent
\textbf{Case A:}
There is $q\in \Delta(A_2)$ such that $u(p,q)\not =\delta x$ for every $p\in \Delta(A_1)$
and every $\delta>0$.

We show that by
playing constantly $q$ (a strategy that we denote by $q^*$) Player 2
can prevent Player~1 from attaining $x$, contradicting the assumption.
Let  $\sigma_1$ be any
strategy of Player 1. Denote by $p_t$ the average mixed action played by
Player 1 up to time $t$,
that is,
$p^t = \frac{1}{t}\int_0^t \sigma_1(s) \rmd s$.
Then, $\gamma^t(\sigma_1,
q^*)=tu(p_t,q)$. Thus, $\gamma^t(\sigma_1, q^*)$ is in the cone
generated by $R_1(q):=\{u(p,q);~p\in \Delta(A_1)\}.$ This cone is closed and
by assumption it does not contain $x$. Thus, there is a positive
distance between $x$ and this cone, implying that
$\gamma^t(\sigma_1, q^*)$ cannot get arbitrarily close to $x$. This
contradicts the fact that the vector $x$ is attainable.

\bigskip
\noindent
\textbf{Case B:} For every $q\in \Delta(A_2)$ there is $p\in \Delta(A_1)$ and $\delta > 0$
such that $u(p,q) =\delta x$, but the $\delta$'s are \emph{not}
bounded away from zero.

In this case, for every $\delta>0$, there is
$q_{\delta}\in \Delta(A_2)$ such that $\delta\ge \max \{\delta';~
\exists p  \hbox{ such that } u(p,q_\delta)=\delta'x\}$. We show that for every
$\delta>0$, if Player 2 plays constantly $q_{\delta}$ (a
strategy that we denote by $q_{\delta}^*$), then there is $\ep>0$
such that for every $\sigma_1$,
$\|\gamma^T(\sigma_1,q_{\delta}^*)-x\|<\ep$ implies
$T>\frac{1}{4\delta}$.

Fix $\delta>0$.
Denote
\[ \delta_0 := \max\{\delta' \colon \exists p \hbox{ such that } u(p,q_\delta) = \delta'x\} < \delta. \]
In particular, $\delta_0x \in R_1(q_\delta)$,
and $\delta'x \not\in R_1(q_\delta)$ for every $\delta' > \delta_0$.
Let $E:=\conv\left(R_1(q_{\delta}) \cup \{\vec
0\} \right)$
be the convex hull of  $R_1(q_{\delta})$
and $\vec 0$.
The set $E$ is convex, compact and it does not contain $\delta' x$ for every $\delta' > \delta_0$.
In particular, $2\delta_0 x\not \in E$. Thus,
there is an open ball $F = B(2\delta_0 x,\eta)$ which is
disjoint of $E$. By the hyperplane separation theorem there is a nonzero
vector $\alpha\in\dR^m$ such that $\langle e, \alpha\rangle \le \langle
f, \alpha \rangle $ for every $e\in E$ and $f\in F$.
Since $\vec 0 \in E$, it follows
that $0= \langle \vec 0 ,\alpha \rangle\le \langle
f,\alpha\rangle$ for every $f\in F$.

Without loss of generality assume that $\|\alpha \|=1$.
We claim that $0 < \langle  x, \alpha\rangle$. Indeed, if  $0 =
\langle x, \alpha \rangle$, then every $f\in F$ can be expressed as
$f=2\delta_0 x+v$, where $v = v(f) \in B(\vec 0,\eta)$.
In particular, $0\le \langle f, \alpha \rangle= \langle
v, \alpha \rangle$.
It follows that $\langle v,\alpha\rangle = 0$ for every $v \in B(\vec 0,\eta)$,
which implies that $\alpha = 0$, contradicting the fact that $\|\alpha\|=1$.

Suppose that $e \in R_1(q_{\delta})$ and $T\cdot e\in B(x,\ep)$,
with $\ep=\langle x,\alpha \rangle /2$. Then, $T\cdot e=x+z$,
where $\|z\|\le \ep$. Thus, $\langle T\cdot e,\alpha
\rangle=\langle x+z,\alpha \rangle $. Since $e\in E$ and
$2\delta_0 x\in F$,
\[ \langle e, \alpha\rangle \le \langle 2\delta_0 x,\alpha\rangle \leq \langle 2\delta
x, \alpha \rangle. \]
Hence,
\begin{equation}\label{eq 111}
T=\frac{\langle x+z,\alpha \rangle }{\langle e,\alpha \rangle}\ge
\frac{\langle x,\alpha \rangle + \langle z,\alpha \rangle}{2\langle \delta x,\alpha \rangle}
\ge \frac{\langle x,\alpha \rangle - \ep}{2\langle \delta x,\alpha \rangle}=\frac{1}{4\delta}.
\end{equation}

Recall that $q_{\delta}^*$ is the strategy of Player 2 that constantly
plays $q_{\delta}$.
To derive a contradiction we will show that the vector $x$ is not attainable; that is, for
every $T$ there is $\ep > 0$ such that for every strategy $\sigma_1$
of Player 1 there is a strategy $\sigma_2$ of Player 2 and $t \leq
T$ satisfying $d(\gamma^t(\sigma_1,\sigma_2),x) > \ep$. Fix a
strategy $\sigma_1$ of Player 1, and suppose that the cumulative
payoff up to time $T$ is within $\ep$ from $x$, that is,
$\|\gamma^T(\sigma_1,q_{\delta}^*)-x\|\le \ep$. Let $p^T :=
\frac{1}{T}\int_0^T \sigma_1(s)\rmd s$ be the average mixed action
played by $\sigma_1$ until time $T$. Thus, $Tu(p
_T,q_{\delta})=x+z$, where $\|z\|\le \ep$. Letting $e=
u(p_T,q_{\delta})$ we obtain by Eq.~(\ref{eq 111}) that
$T>\frac{1}{4\delta}$. In words, the time it takes to reach $B(x,
\ep)$ is at least $\frac{1}{4\delta}$. This shows that there is no
uniform bound on the time at which the total payoff gets close to
$x$. Thus, $x$ is not attainable, which contradicts the assumption.

\bigskip

\noindent\underline{Part 4:}
If Condition \textbf{D4} and Condition \textbf{D1} are satisfied, then Condition \textbf{D3} is satisfied and $x$ is attainable.

We will show that $v_\lambda \geq \delta_0 \langle x,\lambda\rangle$ for every $\lambda \in \dR^m$.
If $\langle x,\lambda\rangle \leq 0$, then by Condition \textbf{D1}
\[
v_{\lambda} \geq 0 \geq \delta_0 \langle x,\lambda\rangle,
\]
as required.
If  $\langle x,\lambda\rangle > 0$ then Condition \textbf{D3} implies that
\[
v_{\lambda}=\inf_{q \in \Delta(A_2)} \sup_{p \in \Delta(A_1)} \langle u(p,q),\lambda\rangle \geq \langle \delta_0 x, \lambda\rangle =
\delta_0 \langle x,\lambda\rangle,
\]
and the proof is complete.

\subsection{Proof of Proposition \ref{suffvec}}

To prove that Conditions $\textbf{E1}$ and $\textbf{E2}$ are sufficient conditions, we prove that they imply Conditions $\textbf{D1}$ and $\textbf{D3}$.
By Corollary \ref{coro2}, Condition \textbf{E1} implies Condition \textbf{D1}.
We now show that Condition \textbf{D3} holds as well.

It is sufficient to prove that Condition \textbf{D3} holds for every $\lambda$ in the unit ball.
The set $\calS^{\geq} := \{\lambda \in \dR^m \colon \|\lambda\|=1, \langle x,\lambda\rangle \geq  0\}$ is compact.
Since the function $\lambda \rightarrow v_\lambda$ is continuous,
Condition \textbf{E2} implies that there exists $\epsilon>0$ such that $v_\lambda \geq \epsilon$
for every
$\lambda \in \calS^{\geq}$.
Let $\delta>0$ such that $\delta \|x\| < \epsilon$.
By Cauchy--Schwartz inequality,
\[
v_\lambda \geq \langle \delta x,\lambda\rangle = \delta \langle x,\lambda\rangle, \ \ \ \forall \lambda \in \calS^{\geq}.
\]

If  $\langle\lambda,x\rangle < 0$ then Condition \textbf{E1} implies that
\[
v_\lambda \geq 0 \geq \delta\langle x,\lambda\rangle,
\]
and the proof is complete.

\end{document}